\documentclass[11pt,letterpaper]{amsart}

 \usepackage{epic,eepic,latexsym, amssymb, amscd, amsfonts, xypic, floatflt}

\input xy
\xyoption{all}

 \newlength{\baseunit}               
 \newcount{\numlines}                
\newcommand{\SYM}{\operatornamewithlimits{\mathsf{sym}}}

\newcommand{\f}{\frac}
\newcommand{\Om}{{\Omega}}
\newcommand{\ld}{{\ldots}}
\newcommand{\bfp}{{\bf p}}
\newcommand{\bx}{{\bf x}}
\newcommand{\bfy}{{\bf y}}
\newcommand{\pa}{{\partial}}
\newcommand{\raw}{{\rangle_g}}
\newcommand{\law}{{\langle}}
\newcommand{\idx}{{\f{x_i\pa}{\pa x_i}}}
\newcommand{\jdx}{{\f{x_j\pa}{\pa x_j}}}
\newcommand{\odx}{{\f{x_1\pa}{\pa x_1}}}

\newcommand{\nodx}{{\f{x_{n+1}\pa}{\pa x_{n+1}}}}
\newcommand{\idy}{\operatorname{\Delta^{\mathit{y}}_{\mathit{i}}}}

\newcommand{\ody}{\operatorname{\Delta^{\mathit{y}}_{1}}}
\newcommand{\nody}{\operatorname{\Delta^{\mathit{y}}_{\mathit{n}+1}}}

\newtheorem{tm}{Theorem}[section]

\newtheorem{co}[tm]{Corollary}

\newcommand{\proj}{\mathbb P}

\newcommand{\sC}{{\mathcal{C}}}
\newcommand{\sF}{{\mathcal{F}}}

\newcommand{\sR}{{\mathcal{R}}}
\newcommand{\sS}{{\mathcal{S}}}
\newcommand{\sT}{{\mathcal{T}}}

\newcommand{\cm}{{\mathcal{M}}}
\newcommand{\cmbar}{\overline{\cm}}

\newcommand{\al}{\alpha}

\newcommand{\be}{\beta}

\newcommand{\si}{\sigma}
\newcommand{\la}{\lambda}
\newcommand{\La}{\Lambda}

\newcommand{\Aut}{\operatorname{Aut}}

\newcommand{\cited}{}

\newcommand{\remind}[1]{{\bf[#1]}}
\newcommand{\lremind}[1]{{\bf[label:  #1]}}
\newcommand{\notation}[1]{}
\renewcommand{\remind}[1]{{}}  
\renewcommand{\lremind}[1]{{}}

\newcommand{\secretnote}[1]{}

\begin{document}
\pagestyle{plain} 
\title{A short proof of the $\lambda_g$-Conjecture without Gromov-Witten
theory: Hurwitz
theory and the moduli of curves}

\author{I. P. Goulden,
D. M. Jackson
and
R. Vakil}
\address{Department of Combinatorics and
Optimization, University of Waterloo
}
\address{Department of Combinatorics and
Optimization, University of Waterloo
}
\address{Department of Mathematics, Stanford University
}

\thanks{The first two authors are partially supported by
NSERC grants.
The third author is partially supported by NSF PECASE/CAREER grant DMS--0238532.
\newline \indent
2000 Mathematics Subject Classification:  Primary 14H10, Secondary 05E99. 
}
\date{April 12, 2006.}
\begin{abstract}
We give a short and direct proof of the $\la_g$-Conjecture.  
The approach is  through the Ekedahl-Lando-Shapiro-Vainshtein theorem, which establishes 
the ``polynomiality'' of Hurwitz numbers, from which we pick off the lowest degree terms. 
The proof is independent of Gromov-Witten theory. 

We briefly describe the philosophy behind our general approach to intersection numbers 
and how it may be extended to other intersection number conjectures. 
\end{abstract}

\maketitle

{\parskip=4pt 

\section{Introduction}

\subsection{Background}
The $\la_g$-Conjecture, now a theorem,  states that
\begin{tm}[The $\lambda_g$-Conjecture]\label{T:lamdeg}
For $n,g\ge 1$,
$$ 
 \int_{\cmbar_{g,n}} \psi_1^{b_1} \cdots \psi_n^{b_n} \la_g  = \binom{2g-3+n} {b_1, \dots, b_n} c_g,$$
where  $\sum_{i=1}^n b_i =2g-3+n$, $b_1,\ld ,b_n\ge 0$ and $c_g$ is a constant that depends only on $g.$
\end{tm}

As usual, 
$\cmbar_{g,n}$ is the (compact) moduli space of stable $n$-pointed genus $g$  curves, 
$\psi_1$, $\dots,$ $\psi_n$ are (complex) codimension $1$ classes corresponding
to the $n$ marked points, and $\la_k$ is the (complex codimension $k$) $k$th Chern
class of the Hodge bundle.  The constant~$c_g$ can be obtained from the $n=1$ case, giving $c_g=
\int_{\cmbar_{g,1}}\psi_1^{2g-2} \la_g$ $= \law \tau_{2g-2} \la_g \raw$, and
throughout the paper $c_g$ is used to denote this particular value. 
For a summary of necessary facts about the moduli space of 
curves, the reader is referred to~\cite{cime}. 
We shall assume background 
about~$\cmbar_{g,n}$ in the Introduction, but the proof of the $\la_g$-Conjecture that is 
presented does not require any knowledge of these notions.

The $\la_g$-Conjecture can be interpreted as a description of the top intersections in the 
tautological cohomology ring of the moduli space ~$\cm_{g,n}^{c}$ of curves of compact type  
(curves whose Jacobian is compact, or equivalently, whose dual graph is a tree).  As such, 
it is part of a family of four problems.  Pandharipande has outlined a philosophy that we should 
expect the ``tautological cohomology rings'' of various moduli spaces to satisfy a ``Gorenstein'' 
property, \emph{i.e.}\ that the top degree term of the ring is one-dimensional, and that the 
multiplication map into it should be a perfect pairing, see~\cite[\S 1]{icm}.  Three spaces 
mentioned there are the moduli space of stable curves~$\cm_{g,n}$, $\cm_{g,n}^c$, and
the moduli space of smooth curves~$\cm_g$ (or, better, the moduli space of pointed curves 
with ``rational tails''~$\cm_{g,n}^{rt}$).  In each case, the one-dimensionality is known 
(see~\cite{socle, fabpjems, thmstar}, for example).

The top intersections in this ring are determined in each case by top intersections 
of~$\psi$-classes  by work of Faber (based on earlier work of Mumford). Then, parallel to 
Pandharipande's  Gorenstein predictions, there are ``intersection-number'' predictions 
determining the full ring  structure.  These are the following: 
\textbf{i)}~the case of $\cmbar_{g,n}$ is \emph{Witten's Conjecture}
(Kontsevich's theorem), which now has a number of very different
and very enlightening proofs;
\textbf{ii)}~the case of $\cm_{g,n}^c$ is the $\la_g$\emph{-Conjecture};  
\textbf{iii)}~the case of $\cm_g$ (or $\cm_{g,n}^{rt}$) is \emph{Faber's intersection number conjecture}.
To these we add a fourth case that seems to be of the same flavor:
\textbf{iv)}~the case of a conjectural compactified universal Picard variety over $\cmbar_{g,n}$ 
(related to double Hurwitz numbers, described in \cite{double}) yields a generating series with
similar behavior (see \cite{double, sz}), which we shall discuss more in~Section~\ref{drphil}.

Our proof of the~$\la_g$-Conjecture is through the Ekedahl-Lando-Shapiro-Vainshtein 
formula~\cite{elsv2},  that establishes the ``polynomiality'' of the Hurwitz numbers, and
by identifying the Hodge integral in the~$\la_g$-Conjecture as a coefficient in the
lowest degree terms in this polynomial.  The proof is short, direct and requires no
Gromov-Witten theory. There are already several proofs of  the 
$\la_g$-Conjecture,  and these will be discussed in 
Section~\ref{lagthm}. 

Our method of proof can be extended to give a proof of Faber's intersection number conjecture
(for up to~3 points, \cite{gjvfaber}). Comments on the philosophy behind this are made in Section~\ref{poga}.

\subsection{Preliminaries}
\subsubsection{The Join-cut Equation}
The Hurwitz numbers~$H^g_\al$  count connected, branched covers  of~$\proj^1$ by a non-singular 
genus $g$ curve,  with branching over~$\infty \in \proj^1$ corresponding to a partition~$\al \vdash d$ 
(these branch  points are ordered), and with simple branching~$(1^{d-2}\,2)$ above~$r=d+n+2g-2$ 
other points,  where~$n=l(\al)$, the number of parts in~$\al$.
Hurwitz~\cite{hur} observed that~$d!H^g_{\al}$ counts the number of factorizations of
an arbitrary permutation in the conjugacy class $\sC_{\al}$ of $\mathfrak{S}_d$ with
cycles of lengths $\al_1, \ld , \al_n$, 
into an ordered,  transitive product of $r$~transpositions in $\mathfrak{S}_d$
(such a product is {\em transitive} if the group generated by the factors
acts transitively on $\{ 1,\ld ,d\}$).

Ordered factorizations are amenable, in principle,  to combinatorial techniques.
The action of a transposition on the disjoint cycles of a permutation can be analyzed by observing that either the transposition joins an~$i$-cycles and 
a~$j$-cycle to make an~$(i+j)$-cycle, or it cuts an~$(i+j)$-cycle into an~$i$-cycle and a~$j$-cycle. 
In this join-cut process, an~$i$-cycle is annihilated by the operator $i\partial/\partial p_i$
and is created by the operator~$p_i$ (regarded as pre-multiplication by~$p_i$) acting on
the~\emph{genus series} 
$$H=\sum_{g\ge 0, n\ge 1} H_n^g x^g,$$
where $H_n^g$ is the \emph{Hurwitz series}, given by
\begin{equation*}
H_n^g(z,\bfp ) =\sum_{d\ge 1}\sum_{{\al\vdash d,} \atop {l(\al )=n}}
|\sC_{\al}|\f{H_{\al}^g}{r!}p_{\al}z^d,
\end{equation*}
with $\al=(\al_1, \ld ,\al_n)$ and $p_{\al}=p_{\al_1}\cdots p_{\al_n}$. 
It follows immediately from this construction that the genus series
satisfies the \emph{Join-cut Equation}~(see~\cite{gjvai}):
\begin{eqnarray}\label{e:jc}
\biggl( z\f{\pa}{\pa z} +2x\f{\pa}{\pa x}-2 &+& \sum_{i\ge 1}p_i\f{\pa}{\pa p_i} 
 \biggr) H\\
&=&\tfrac{1}{2}\sum_{i,j\ge 1}
\left(ij x p_{i+j}\f{\pa^2H}{\pa p_i\pa p_j}+ijp_{i+j}\f{\pa H}{\pa p_i}\f{\pa H}{\pa p_j}\right.
+\left. (i+j)p_ip_j\f{\pa H}{\pa p_{i+j}}\right) ,\notag
\end{eqnarray}
where the first two operators on the right hand side give
the cycle-type after a~\emph{join} and the third
operator gives the cycle-type after a~\emph{cut}. Because of transitivity, there are two
cases of joins.
The first operator is a  join of two cycles within a \emph{single} transitive factorization, 
while the second operator is
a join of two cycles, one from each of \emph{two} disjoint transitive ordered factorizations.

\subsubsection{The Genus Expansion Ansatz}
The background to our proof is an observation about Hurwitz 
numbers~$H^g_{\al}.$  For fixed~$n=l(\al)$ and~$g$, with~$n,g\ge 1$ or~$n\ge 3,g=0$, 
it was conjectured that
\begin{equation}\label{cosypoly}
H^g_{\al} = r!  \prod_{i=1}^n \left(  \frac { \al_i^{\al_i}} {\al_i!} \right) P_{g,n}(\al_1, \dots, \al_n),
\end{equation}
for some symmetric polynomial~$P_{g,n}$ in the~$\al_i$, with terms of total degrees 
between~$2g-3+n$ and~$3g-3+n$. This important property is essentially the \emph{Polynomiality
Conjecture} of~\cite[Conj.\ 1.2]{gjconj} (the connection is 
made in~\cite{gjv1}).
The Polynomiality Conjecture was settled by Ekedahl, Lando, M. Shapiro, and Vainshtein, 
who proved the remarkable ELSV-formula \cite{elsv1, elsv2}. (For a proof in the context 
of Gromov-Witten theory, see \cite{gvelsv}, and also \cite{thmstar}.)
In the present notation, the ELSV-formula states that\lremind{elsvf}
\begin{equation}\label{elsvf}
P_{g,n} = \int_{\cmbar_{g,n}} \frac {1 - \la_1 + \cdots + (-1)^g \la_g}
{ (1 - \al_1 \psi_1) \cdots (1 - \al_n \psi_n)}.
\end{equation}
Equation~\eqref{elsvf} should be interpreted as follows:
formally invert the denominator as a geometric series;
select the terms of codimension~$\dim \cmbar_{g,n}=3g-3+n$;
and ``intersect'' these terms on $\cmbar_{g,n}$.
The formula therefore yields
\begin{equation}\label{geoexp}
P_{g,n} = \sum_{{b_1 + \cdots + b_n + k = 3g-3+n,} \atop {b_i \geq 0, \; 0\leq k\leq g}} 
(-1)^k 
\law\tau_{b_1} \cdots \tau_{b_n}\la_k\raw
\al_1^{b_1} \cdots \al_n^{b_n},
\end{equation}
where we have used the Witten symbol
(from Gromov-Witten theory)
$$\law\tau_{b_1} \cdots \tau_{b_n}\la_k\raw :=
\int_{\cmbar_{g,n}} \psi_1^{b_1} \cdots \psi_n^{b_n} \la_k,$$
and note that
\begin{equation}\label{witzero}
\law\tau_{b_1} \cdots \tau_{b_n}\la_k\raw =0
\end{equation}
unless $b_1+\cdots+b_n=3g-3+n-k$.

Then substituting ~\eqref{geoexp} into~\eqref{cosypoly}, we
obtain the {\em Genus Expansion Ansatz} for the 
Hurwitz series  (see Thm.~2.5 of~\cite{gjv1} for details), namely
\begin{equation}\label{hurwwit}
H^g_n=\f{1}{n!}
\sum_{{b_1, \ld  ,b_n\ge 0,} \atop {0\le k\le g}}
(-1)^k\law\tau_{b_1} \cdots \tau_{b_n}\la_k\raw
\prod_{i=1}^n\phi_{b_i}(z,\bfp )
\end{equation}
for $g,n \ge 1$ and for $g=0$, $n\ge 3$, where
\begin{equation*}
\phi_i(z,\bfp )=\sum_{m\ge 1}\f{m^{m+i}}{m!}p_mz^m,\quad i\ge 0.
\end{equation*}
This should be interpreted as just a re-writing of the ELSV formula.

\subsubsection{Our approach to the~$\la_g$-Conjecture}
The second observation about  $P_{g,n}(\al)$ (recall that the first is that it is a polynomial)
is that its lowest total degree (this is~$2g-3+n$) part appears to have the form
$$(\al_1 + \cdots + \al_n)^{2g-3+n} c_g,$$ 
where $c_g$ is a constant depending only upon~$g.$ 
This assertion is equivalent to the~$\la_g$-Conjecture
by \eqref{geoexp}
and is the form of the result that we prove.

We require only two properties of the Hurwitz series, 
namely that it satisfies the Join-cut Equation~\eqref{e:jc} and that it has the Genus Expansion 
Ansatz~\eqref{hurwwit}. 
To obtain a characterization of the left hand side of Theorem~\ref{T:lamdeg} in terms of an operator
acting on the Hurwitz series $H^g_n$
we transform the latter in a series of three steps:
\begin{enumerate}
\item[\textbf{(i)}]   symmetrization of the Hurwitz  series and the Join-cut Equation; 
\item[\textbf{(ii)}]   change of variables to obtain a polynomial; and 
\item[\textbf{(iii)}]   determination of the full (to be defined later) terms
 of minimum degree in this polynomial. 
\end{enumerate}

In Section~\ref{s:gea}, we apply this transformation to
the Genus Expansion Ansatz for the Hurwitz series. In our main result
of this section, Theorem~\ref{minhur}, we prove that each Witten
symbol whose evaluation is the subject of the $\la_g$-Conjecture is
the coefficient of a unique monomial in the transformed Hurwitz series.
In Section ~\ref{s:jc}, we apply this transformation to the Join-cut
Equation \eqref{e:jc} for the Hurwitz series.  In our 
main result of this
section, Theorem~\ref{newpde}, we prove that a genus generating series
for the transformed Hurwitz series satisfies a simple
partial differential equation. We then solve this partial differential
equation in Theorem~\ref{solpde}.  Finally, in Section ~\ref{s:pf},
we prove the $\la_g$-Conjecture by comparing the results
obtained in Sections ~\ref{s:gea} and ~\ref{s:jc}.

We note in passing that the transformations we apply in this paper are also used
in~\cite{gjvfaber} in which we are able
to prove (up to $3$ parts) the Faber intersection number conjecture (see~\cite{fabconj}). In the
latter Faber case, we apply the steps~\textbf{(i)} and~\textbf{(ii)} of the transformations
applied in the present paper, but for step~\textbf{(iii)}, in the Faber case, we 
consider terms of {\em maximum} degree rather than the minimum degree
(on a {\em different} polynomial).
This philosophy will be discussed
in Section~\ref{poga}.

In the Appendix, we indicate how our approach can be used to obtain the generating series
of intersection numbers that are close to ``minimum'' in the sense that has been described above,
and we exhibit the explicit series in a few cases.

\subsection{Previous proofs of the $\la_g$-Conjecture}
\label{lagthm}\lremind{lagthm}
The $\la_g$-Conjecture was first proved in Faber and Pandharipande's landmark 
paper~\cite{fabplambdag}.  Their approach was to use localization on the space of stable maps 
to~$\proj^1$ to obtain relations among these intersection numbers.  They then showed 
that the~$\la_g$-Conjecture's prediction satisfied these relations. Finally, they proved that the 
relations uniquely determined the predictions of the $\la_g$-Conjecture by establishing
the invertibility of  a large matrix whose entries are counts of various partitions; this requires 
seven pages of explicit calculation.

A second proof is as follows.  Getzler and Pandharipande showed that
the $\la_g$-Conjecture is a formal consequence of the Virasoro
Conjecture for $\proj^1$~\cite[Thm.~ 3]{gp}, by showing that it satisfies a recursion arising from the
Virasoro Conjecture, and then showing that the recursion has a unique
solution.  The Virasoro Conjecture for~$\proj^1$ was then
shown in two ways.  It was
proved for all curves by Okounkov and Pandharipande~\cite{op}.  Also,
Givental has announced a proof of the Virasoro Conjecture for Fano
toric varieties~\cite{giv}.  The details have not yet appeared, but
Y.-P. Lee and Pandharipande are writing a book \cite{leep} supplying
them.  These proofs of the Virasoro Conjecture in important cases are
among the most significant results in Gromov-Witten theory, and this
method of proof of the~$\la_g$-Conjecture seems somewhat circuitous.
(Much of this paragraph also applies to Faber's intersection number
conjecture.)

Liu, Liu, and Zhou gave a new proof in~\cite{llz} as a consequence of the Mari\~no-Vafa 
formula~\cite{mv}, which was proposed by the physicists Mari\~no and Vafa
and proved by  Liu, Liu, and Zhou in \cite{llz1}.  
This Gromov-Witten-theoretic proof is quite compact.

\section{Transformation of the Genus Expansion Ansatz}
\label{s:gea}\lremind{s:gea}
In this section, we transform the Hurwitz series $H^g_n$ through the Genus Expansion 
Ansatz~\eqref{hurwwit} by constructing the operator to extract the intersection number of 
Theorem~\ref{T:lamdeg}. 

\newcommand{\Xin}{\operatorname{\Xi_\mathit{n}}}
\renewcommand{\sC}{\operatorname{\mathsf{C}}}
\renewcommand{\sF}{\operatorname{\mathsf{F}}}

\newcommand{\sFm}{\operatorname{\mathsf{M}}}
\newcommand{\oO}{\operatorname{\mathsf{O}}}
\newcommand{\oE}{\operatorname{\mathsf{E}}}

\subsection{Step~1 -- Symmetrization}
For the first step of our transformation, we {\em symmetrize} the Hurwitz series using
the linear symmetrization operator $\Xin$, given by
\begin{equation*}
\Xin p_{\al}z^{|\al|}=\sum_{\si\in \mathfrak{S}_n}x_{\si(1)}^{\al_1} \cdots x_{\si(n)}^{\al_n},
\quad n\ge 1,
\end{equation*}
if $l(\al )=n$ (with $\al=(\al_1,  \ld  ,\al_n)$), and $0$ otherwise.
Thus, applying $\Xin$ to~\eqref{hurwwit} we obtain,
for $n,g\ge 1$ and $n\ge 3,g\ge 0$,
\begin{equation}\label{Psimg}
\Xin H_n^{g}=\f{1}{n!}
\sum_{{b_1, \ld  ,b_n\ge 0,} \atop {0\le k\le g}}
(-1)^k\law\tau_{b_1} \cdots \tau_{b_n}\la_k\raw
\sum_{\si\in \mathfrak{S}_n}\prod_{i=1}^n\phi_{b_i}(x_{\si(i)}),
\end{equation}
where 
$$\phi_i(x)=\phi_i(x,{\bf 1})=\sum_{m\ge 1}\f{m^{m+i}}{m!}x^m.$$
We note that  
\begin{equation}\label{phiw}
\phi_i(x)=\left(x\f{d}{dx}\right)^{i+1}w(x),
\end{equation}
where 
$$w(x)=\sum_{m\ge 1} m^{m-1}\f{x^m}{m!}$$ 
is the (exponential) generating series for the number~$m^{m-1}$ of trees with~$m$ vertices,
labelled from~$1$ to~$m,$ and having a single vertex which is further distinguished (for example, by
painting it red).
Such trees are termed vertex-labelled rooted trees, and we shall refer to~$w(x)$ as the
\emph{rooted tree series}.  It is the unique formal power series solution of the (transcendental)
functional equation~(see \emph{e.g.}~\cite{gjce}~\S3.3.10)
\begin{equation}\label{rotreq}
w=x e^{w}
\end{equation}
(which we shall refer to as the {\em rooted tree equation}).

\newcommand{\spn}{\mathsf{span}}

\subsection{Step~2 -- change of variables}
We next consider a change of variables for the symmetrized Hurwitz series.
Consider $y(x)=(1-w(x))^{-1}$. Then 
\begin{equation}\label{compinv}
y(x)=1+\sum_{m\ge 1}\f{m^m}{m!}x^m= 1+\phi_0(x),
\end{equation}
which can be seen most easily perhaps from~(\ref{opxw}) below.
Let $w_j=w(x_j)$ and $y_j=y(x_j)$, $j=1,\ld ,n$, and
let $\sC$ be an operator, applied to a formal power series in $x_1,\ld ,x_n$,
that changes variables, from  the indeterminates $x_1,\ld ,x_n$ to $y_1,\ld ,y_n$.
Thus, from~(\ref{compinv}), to carry out $\sC$ we substitute $x_j=g(y_j-1)$, where $g$ is the
compositional inverse of $\phi_0$.
In general, this will not yield a formal power series in $y_1,\ld ,y_n$, but
when we apply $\sC$ to $\Xin H^g_n$, we do obtain a formal power series (in fact,
for each fixed $n,g$ it is a polynomial) as we prove below.

First we prove some properties of $\sC$.
Differentiating the rooted tree equation~(\ref{rotreq}), we obtain the operator identity
\begin{equation}\label{opxw}
x_j\f{d}{dx_j}=\f{w_j}{1-w_j}\f{d}{dw_j}.
\end{equation}
But $dy_j=y_j^2dw_j$, so we have the operator identities
\begin{equation}\label{commut}
\sC\jdx=(y_j^3-y_j^2)\f{\pa}{\pa y_j}\sC,\quad\quad
\sC w_j\f{\pa}{\pa w_j} = (y_j^2-y_j)\f{\pa}{\pa y_j}\sC,
\end{equation}
where when we apply $\sC$ to expressions involving $w_j$, we
interpret $w_j$ as $w(x_j).$
{}From~(\ref{phiw}),~(\ref{opxw}) and~(\ref{commut}), we also obtain
\begin{equation}\label{phiy}
\sC\phi_i(x_j) 
\left( (y_j^3-y_j^2)\f{\pa}{\pa y_j}\right)^i(y_j-1), \qquad \mbox{for $i\ge 0$.}
\end{equation}

Now~(\ref{witzero}),~(\ref{Psimg}) and~(\ref{phiy}) together enable us to obtain
a polynomial expression for $\sC\,\Xin H_n^g$. The fact that this is unique,
and hence that the application of $\sC$ to $\Xin H^g_n$ is well-defined for
formal power series, follows immediately from the fact that the
non-negative powers of the rooted tree series $w(x)$ are linearly independent,
as formal power series in $x$.

\subsection{Step~3 -- full terms of minimum total degree}
The final step in the transformation of the Hurwitz series is to identify a particular
subset of terms.
We say that a monomial $y_1^{i_1} \cdots y_n^{i_n}$ is {\em full} if $i_1, \ld  ,i_n\ge 1$.  
Let $\sF_k f$ be the subseries of a series $f$ in $y_1, \ld  ,y_n$ consisting of the full terms 
of total degree~$k.$ Thus, for example, from~(\ref{phiy}), we immediately obtain
\begin{equation}\label{Fpsii}
\sF_{i+1}\sC\phi_i(x_j)=(-1)^{i}i!y_j^{i+1},
\end{equation}
by induction on $i\ge 0$, and
\begin{equation}\label{Fpsik}
\sF_{k}\sC\phi_i(x_j)=0, \quad i\ge 0,\;\; k<i+1.
\end{equation}
In addition, when applied to $\sC\,\Xin H_n^g$, let $\sFm$ denote $\sF_{2g-3+2n}.$

Let $m_{\be}$ denote the monomial symmetric function, where we allow $0$ parts in $\be$, 
and write $\be\vdash_0 d$ to indicate that $\be$, with parts equal to $0$ allowed, is a partition 
of $d$.  As usual, $l(\be)$ is the number of parts of $\be$ (including the parts equal to $0$).

\begin{tm} \label{minhur}\lremind{minhur}
Let $\bfy=(y_1,\ldots,y_n)$.
For $n,g\ge 1$  and $n\ge 3, g=0$,
\begin{equation*}
\sFm \sC\,\Xin H_n^g=y_1\cdots y_n(-1)^{3g-3+n}\!\!\!\!
\sum_{{\be\vdash_0 2g-3+n,} \atop {l(\be)=n}}\!\!\!\!
\law\tau_{\be_1} \cdots \tau_{\be_n}\la_g\raw \be_1! \cdots \be_n!m_{\be}(\bfy ),
\end{equation*}
where $\be=(\be_1, \ld  ,\be_n)$, and
\begin{equation}\label{zeromin}
\sF_{k}\sC\,\Xin H_n^g=0,\quad\text{for}\;\; k<2g-3+2n.
\end{equation}
\end{tm}

\begin{proof}
We apply $\sFm \sC$ to the symmetrized Genus Expansion Ansatz~(\ref{Psimg}), so
from~(\ref{witzero}),~(\ref{Fpsii}) and~(\ref{Fpsik}), we obtain
\begin{eqnarray*}
\sFm \sC\,\Xin H^g_n
= \f{1}{n!}\!\!\!\!\!\!
\sum_{{b_1, \ld  ,b_n\ge 0,} \atop {b_1+ \cdots +b_n=2g-3+n}}
\!\!\!\!\!\!\!\!\!\!\!\!
(-1)^g\law\tau_{b_1} \cdots \tau_{b_n}\la_g\raw
\sum_{\si\in \mathfrak{S}_n}\prod_{i=1}^n (-1)^{b_i}b_i!y_{\si(i)}^{b_i+1}. 
\end{eqnarray*}
But we have
\begin{equation*}
\sum_{\si\in \mathfrak{S}_n}\prod_{i=1}^n y_{\si(i)}^{b_i}=
\left| \Aut\be \right| m_{\be}(\bfy ),
\end{equation*}
where $\be$ is the partition (with $0$ allowed as parts)
whose parts are $b_1, \ld  ,b_n$, reordered and $\Aut \be$ is
the subgroup of $\mathfrak{S}_n$
preserving $(b_1, \dots, b_n)$ (through its permutation action on the coordinates).
The first part follows by changing the range of summation from $b_1, \ld  ,b_n$ to $\be$.
The second part  follows immediately from~(\ref{Psimg}) and~(\ref{Fpsik}).
\end{proof}

Note that~(\ref{zeromin}) implies that there are no full terms
in the series $\sC\,\Xin H^g_n$ whose total degree is less than $2g-3+2n$.
Thus we say that $\sFm \sC\,\Xin H^g_n$ consists of the full terms of {\em minimum} total degree
in $\sC\,\Xin H_n^g$ (though we understand that this is informal, since
it assumes that the full terms of total degree $2g-3+2n$ are
not identically zero).

\noindent
{\bf An aside.}  Theorem~5.1 of~\cite{gjvai}, which is not used in this paper, concerns the 
terms of \emph{maximum} total degree when we apply $\sC$ since it gives
 an upper limit for the total degree.
This should be corrected. The total degree of the terms is in fact less than or
equal to $3m-6+3g$, {\em not} $2m-5+6g$, as was incorrectly given there.

\section{Transformation of the Join-cut Equation}
\label{s:jc}\lremind{s:jc}

We now apply the operator~ $\sFm \sC\,\Xin$ to the Join-cut Equation  \eqref{e:jc}
to derive a partial differential equation for the generating series for $\sFm \sC\,\Xin H^g_n$.
The following notation is needed for  applying the symmetrization 
operator $\Xin$ and stating the equation.

For $i,j\ge 0$, $i+j\le
n$, let 
$\underset{i,j}{\overset{x}{\SYM}}$
be the mapping, applied to a series in $x_1,
\ld ,x_n$, given by
\begin{equation*}
\SYM_{i,j}^x f(x_1, \ld  ,x_n)=\sum_{\sR,\sS,\sT}f(\bx_{\sR},\bx_{\sS},\bx_{\sT}),
\end{equation*}
where the sum is over all ordered partitions $(\sR,\sS,\sT)$ of $\{1, \ld  ,n\}$,
where
$\sR=\{x_{r_1}, \ld  ,x_{r_i}\}$, $\sS=\{x_{s_1}, \ld  ,x_{s_j}\}$,
$\sT=\{x_{t_1}, \ld  ,x_{t_{n-i-j}}\}$
and $(\bx_{\sR},\bx_{\sS},\bx_{\sT})=(x_{r_1}, \ld  ,x_{r_i},x_{s_1}, \ld  ,x_{s_j},x_{t_1}, \ld  ,x_{t_{n-i-j}})$,
and where $r_1< \ldots <r_i$, $s_1< \ldots <s_j$, and $t_1< \ldots <t_{n-i-j}.$
If $i$ or $j$ is equal to $0$, then we may suppress them by writing $\underset{2}{\overset{x}{\SYM}}$ for $\underset{2,0}{\overset{x}{\SYM}}$,
for example.

The following result gives an expression for the result of applying the symmetrization 
operator~$\Xin $ to the Join-cut Equation for the Hurwitz series.

\begin{tm}[see~\cite{gjvai}~Thm.~4.4]\label{pdesymx}
The series $\Xin H_n^g$ satisfy the partial differential equation
\begin{equation*}
\left(\sum_{i=1}^nw_i\f{\pa}{\pa w_i} +n+2g-2\right)\Xin H_n^{g}(x_1, \ld  ,x_n)
=T_1+ T_2+T_3 +T_4,
\end{equation*}
where
\begin{eqnarray*}
T_1&=&\tfrac{1}{2}\sum_{i=1}^n
\left.\left(\idx\,\nodx\,\Xin H_{n+1}^{g-1}(x_1, \ld  ,x_{n+1})\right)\right|_{x_{n+1}=x_i},
\\
T_2&=&\SYM_{1,1}^x\f{w_2}{1-w_1}\f{1}{w_1-w_2}\odx\Xin H_{n-1}^{g}(x_1,x_3, \ld  ,x_{n}),
\\
T_3&=&\sum_{k=3}^n\,\SYM_{1,k-1}^x
\left(\odx\Xin H_{k}^{0}(x_1, \ld  ,x_{k})\right)
\left(\odx\Xin H_{n-k+1}^{g}(x_1,x_{k+1}, \ld  ,x_{n})\right),   
\\
T_4&=&\tfrac{1}{2}\sum_{{1\le k\le n,} \atop {1\le a\le g-1}}
\SYM_{1,k-1}^x\left(\odx\Xin H_{k}^{a}(x_1, \ld  ,x_{k})\right)
\left(\odx\Xin H_{n-k+1}^{g-a}(x_1,x_{k+1}, \ld  ,x_{n})\right),
\end{eqnarray*}
for $n,g\ge 1$,
with initial condition $\Xin H_0^g=0$ for $g\ge 1$.
\end{tm}

Here we shall 
only consider Theorem~\ref{pdesymx} for $n\ge 1,g\ge 2$ and $n\ge 2, g=1$,
and note that for this range of values, $\Xin H^0_i$ only arises in
this equation for $i\ge 3$.
In the statement of the result, a meaning is attached to $(w_i-w_j)^{-1}$ 
for $1\le i<j\le n$ by
imposing the total order $w_1\prec\ldots\prec w_n$, and then
defining $(w_i-w_j)^{-1}$ $=-w_j^{-1}(1-w_i/w_j)^{-1}$.
This  then defines a formal power series ring in $w_1$ with coefficients
that are formal Laurent series in $w_2,\ld ,w_{n}$ (see Xin~\cite{x}).

We now consider the partial differential equation for a genus
generating series $\Om_n(y_1, \ld  ,y_n;t)$,  which arises
by applying $\sFm \sC$ to the symmetrized Join-cut Equation
given in Theorem~\ref{pdesymx}. For this purpose,
let $\oO f$ and $\oE f$ denote, respectively, the {\em odd}
and {\em even} subseries of the formal power series $f$ in
the indeterminate $t$.

\begin{tm}\label{newpde}
Let
\begin{equation*}
\Om_n(\bfy ;t)=\sum_{g\ge 1}\f{(-1)^{3g-3+n}}{c_g} \sFm \sC\,\Xin H_n^g
\f{t^{2g-3+n}}{(2g-3+n)!},\quad n\ge 1.
\end{equation*}
Then, for $n\ge 2$, we have the partial differential equation
\begin{equation*}
(n-1)\f{\pa}{\pa t}\Om_n(\bfy ;t)=
\SYM^y_{1,1}\f{y_1^3y_2}{y_1-y_2}\f{\pa}{\pa y_1}\Om_{n-1}(y_1,y_3, \ld  ,y_n;t),
\end{equation*}
with initial condition $\Om_1(y_1;t)=\oE\,\f{y_1}{1-y_1t}$.
\end{tm}
\begin{proof}
We begin by applying $\sC$ to Theorem~\ref{pdesymx}, and note that
\begin{equation*}
\sC\f{w_2}{1-w_1}\f{1}{w_1-w_2}=y_1^2\f{y_2-1}{y_1-y_2}.
\end{equation*}
Let $ \Delta_j^y=(y_j^3-y_j^2)\f{\pa}{\pa y_j}.$
Then this result, together with~(\ref{commut}), transforms the equation in
Theorem~\ref{pdesymx} into
a partial differential equation for $\sC\,\Xin H_n^g$ given by
\begin{equation}\label{pdesymy}
\left(\sum_{i=1}^ny_i(y_i-1)\f{\pa}{\pa y_i} +n+2g-2\right)\sC\,\Xin H_n^{g}(y_1, \ld  ,y_n)
=T'_1+ T'_2 + T'_3+ T'_4,
\end{equation}
where $n\ge 2, g=1$ or $n\ge 1, g\ge 2$, and
\begin{eqnarray*}
T'_1&=&\tfrac{1}{2}\sum_{i=1}^n
\left.\left(\idy\nody\sC\,\Xin H_{n+1}^{g-1}(y_1, \ld  ,y_{n+1})\right)\right|_{y_{n+1}=y_i},
\\
T'_2&=&\SYM^y_{1,1} y_1^2\f{y_2-1}{y_1-y_2}\ody\sC\,\Xin H_{n-1}^{g}(y_1,y_3, \ld  ,y_{n}),
\\
T'_3&=&\sum_{k=3}^n\,\SYM^y_{1,k-1}
\left(\ody\sC\,\Xin H_{k}^{0}(y_1, \ld  ,y_{k})\right)
\left(\ody\sC\,\Xin H_{n-k+1}^{g}(y_1,y_{k+1}, \ld  ,y_{n})\right),
\\
T'_4&=&\tfrac{1}{2}\sum_{{1\le k\le n,} \atop {1\le a\le g-1}}
\SYM^y_{1,k-1}\left(\ody\sC\,\Xin H_{k}^{a}(y_1, \ld  ,y_{k})\right)
\left(\ody\sC\,\Xin H_{n-k+1}^{g-a}(y_1,y_{k+1}, \ld  ,y_{n})\right).
\end{eqnarray*}
Now apply $\sFm $ to~(\ref{pdesymy}), and
use~(\ref{zeromin}). With the notation $\Om_n^g=\sFm \sC\,\Xin H_n^g$, 
the only non-zero contributions on the left hand side arise from
\begin{equation*}
\left(-\sum_{i=1}^ny_i\f{\pa}{\pa y_i}+n+2g-2\right) \Om_n^g=
\left(-(2g-3+2n)+n+2g-2\right)\Om_n^g =(1-n)\Om_n^g,
\end{equation*}
since all terms in $\Om_n^g$ have total degree $2g-3+2n$.
On the right hand side, all contributions from
terms $T'_1$, $T'_3$ and $T'_4$ are zero.
For $T'_2$, the only non-zero contributions arise from
\begin{equation*}
\SYM^y_{1,1} \f{y_1^4}{y_1-y_2}\f{\pa}{\pa y_1}\Om^g_{n-1}(y_1,y_3, \ld  ,y_{n}),
\end{equation*}
from degree considerations alone.
However, note that $y_1^4/(y_1-y_2)=y_1^3+y_1^3y_2/(y_1-y_2)$,
and we conclude that, for {\em full} terms, the non-zero contributions
from $T'_2$ are given by
\begin{equation*}
\SYM^y_{1,1} \f{y_1^3y_2}{y_1-y_2}\f{\pa}{\pa y_1}\Om^g_{n-1}(y_1,y_3, \ld  ,y_{n}).
\end{equation*}
Thus, we obtain the partial differential equation
\begin{equation}\label{dfeqzero}
(1-n)\Om_n^g(\bfy )=
\SYM^y_{1,1} \f{y_1^3y_2}{y_1-y_2}\f{\pa}{\pa y_1}\Om^g_{n-1}(y_1,y_3, \ld  ,y_{n}),
\end{equation}
for $n\ge 2,g\ge 1$.
\medskip

Now multiply this equation for $\Om^g_n$ by $(-1)^{3g-4+n}t^{2g-4+n}/c_g(2g-4+n)!$,  and
sum over $g\ge 1$, to obtain the partial differential equation for $\Om_n$, $n\ge 2$.
For $n=1$, we have 
\begin{equation*}
\Om_1^g= (-1)^{3g-2}\law\tau_{2g-2}\la_g\raw (2g-2)! y_1^{2g-1},
\end{equation*}
from Theorem~\ref{minhur}, which gives $\Om_1(y_1;t)=\sum_{g\ge 1} y_1^{2g-1}t^{2g-2}$, and
the result follows. 
\end{proof}

The partial differential equation in Theorem~\ref{newpde} is simple enough
that it can be solved explicitly.

\begin{tm}
For $n\ge 1$,\label{solpde}\lremind{solpde}
$$
\Om_n(\bfy ;t)= \left\{ \begin{array}{ll}
\oE\, \displaystyle{\prod_{i=1}^n\f{y_i}{1-y_it}}, &  \text{for}\;\; n\;\;\text{odd},\\
\oO\, \displaystyle{\prod_{i=1}^n\f{y_i}{1-y_it}}, &  \text{for}\;\; n\;\;\text{even}.
\end{array} \right.
$$
\end{tm}

\begin{proof}
Let $F_n(\bfy ;t)=\prod_{i=1}^n\f{y_i}{1-y_it}$. Then we have
\begin{eqnarray*}
\SYM^y_{1,1}\f{y_1^3y_2}{y_1-y_2}\f{\pa}{\pa y_1}F_{n-1}(y_1,y_3, \ld  ,y_n;t)
&=&F_n(\bfy ;t)\SYM^y_{1,1}\f{y_1^2(1-y_2t)}{(y_1-y_2)(1-y_1t)}.
\end{eqnarray*}
But the symmetrized term on the right hand side of this equation becomes
\begin{eqnarray*}
\SYM^y_{2}\f{y_1^2(1-y_2t)^2-y_2^2(1-y_1t)^2}{(y_1-y_2)(1-y_1t)(1-y_2t)} 
= \SYM^y_{2}\left(\f{y_1}{1-y_1t}+\f{y_2}{1-y_2t}\right) 
= (n-1)\sum_{i=1}^n\f{y_1}{1-y_1t},
\end{eqnarray*}
and we thus have
\begin{equation*}
\SYM^y_{1,1}\f{y_1^3y_2}{y_1-y_2}\f{\pa}{\pa y_1}F_{n-1}(y_1,y_3, \ld  ,y_n;t)=
(n-1)\f{\pa}{\pa t}F_n(\bfy ;t).
\end{equation*}
This proves that $F_n(\bfy ;t)$ is a solution to the partial
differential equation given in Theorem~\ref{newpde}, and the result
follows from the initial conditions and the parity restrictions
on the generating series $\Om_n(\bfy ;t)$.
\end{proof}

\section{Proof of the $\la_g$-Conjecture}
\label{s:pf}\lremind{s:pf}
Now we can prove the $\la_g$-Conjecture stated as~Theorem~\ref{T:lamdeg}.

\begin{proof}
We have $\prod_{i=1}^n(1-y_it)^{-1}=\sum_{k\ge 0}h_k(\bfy )t^k$,
where $h_k(\bfy )$ is the $k$th {\em complete} (or {\em homogeneous})
symmetric function, given by 
\begin{equation*}
h_k(\bfy )=\sum_{{\al\vdash_0 k,} \atop {l(\al)=n}}m_{\al}(\bfy ).
\end{equation*}
Then, immediately from Theorem~\ref{solpde}, we obtain
\begin{equation*}
\Om^g_n(\bfy )=c_g(-1)^{3g-3+n}(2g-3+n)!\;\; y_1 \cdots y_n\!\!\!\!
\sum_{{\al\vdash_0 2g-3+n,} \atop {l(\al)=n}}m_{\al}(\bfy ).
\end{equation*}
and the result follows by comparing this result with Theorem~\ref{minhur}.
\end{proof}

\section{The philosophy of the general approach}\label{poga}
The approach stands in a more general geometric-combinatorial setting, and although we do not need much of this setting here, we do require it for our proof 
\cite{gjvfaber} of Faber's intersection number conjecture (for a small number of points). This more general setting provides a useful perspective for the proof that we have given of the $\la_g$-Conjecture.

\subsection{A bridge between geometry and combinatorics}
The general approach is based on the observation that localization theory
(developed in Gromov-Witten theory by \cite{vl}), when applied to
the cases that have been described above, expresses a series in the intersection
numbers in terms of a sum over combinatorial structures (such as trees or graphs) that
are weighted by Hurwitz numbers~$H_\al^g$  (or double Hurwitz numbers in the
case of Faber's Conjecture).  An account of this is given in~\cite{cime}.
In this sense,  localization theory provides a bridge from the \emph{geometry} of 
intersection numbers for the moduli
spaces of curves on the one hand, to branched covers on the other. As we have seen, the latter
may be regarded as \emph{combinatorial} structures. 

Associated with the generating series for transitive ordered factorizations into transpositions
is a functional equation that leads to an implicitly defined set of series.  These, together
with the combinatorial structure (trees, graphs) that are a consequence of the use of localization theory, determine 
an implicit change of variables.  Although the functional equation is transcendental, the derivatives
of its solution are, in effect,  \emph{rational} in the solution.  It is precisely this rationality that leads to the 
\emph{polynomiality property} and thence to a linear system of equations for the intersection numbers.

The usefulness of this general point of view is reinforced by the
following observations.  First, it enables us to obtain other Hodge
integrals. Secondly, our proof of Faber's intersection number
conjecture (for a small number of points) uses localization theory to
create a sum over a particular class of trees weighted by genus $0$
double Hurwitz numbers, which we subject
to a similar but more complex (combinatorial) analysis.

\subsection{Integrable systems, recent developments and closing comments}
\label{drphil}\lremind{drphil}

The $\la_g$-Conjecture, a statement about the moduli space of curves, or the factorization of 
transpositions, should not need to follow from Gromov-Witten theory. 
This work was motivated by the 
fact that the other three intersection-number conjectures either follow or might be expected to 
follow from understanding the algebraic structure of Hurwitz-type numbers.  In each case, there 
is a natural change of variables (motivated by the string and dilaton equations); and in each case, 
there is a connection to integrable hierarchies. We point out the following recent developments:
\textbf{i)} Kazarian and Lando's~\cite{kal} and Kim and Liu's~\cite{kil} 
short proofs
of Witten's conjecture (the $\cmbar_{g,n}$ case);
\textbf{ii)} Shadrin and Zvonkine's description and proof of a Witten-type theorem
on the conjectural compactified Picard variety (related to one-part double Hurwitz numbers), 
relating the intersection theory to integrable hierarchies~\cite{sz}; and
\textbf{iii)} our proof of Faber's intersection number conjecture for up to three points, using 
``Faber-Hurwitz numbers,'' \cite{gjvfaber}.

Finally, the Join-cut Equation seems intertwined in
some way with integrable hierarchies, but the precise connection is
not yet clear.  For example, it is a non-trivial task to go from the
Join-cut Equation to Witten's Conjecture.
\vspace{.1in}


\noindent {\large{\bf Acknowledgments.}}
We thank R. Cavalieri, S. Lando, S. Shadrin and
D. Zvonkine for comments which have improved the manuscript,
and the third author thanks T. Graber for helpful conversations.

\appendix

\section{Intersection numbers $k$ higher than minimum}\label{S:oin}
In principle, the formalism that we have described can be used also to obtain
$\law\tau_{\al_1} \cdots \tau_{\al_n}\la_{g-k}\raw$ for $k>0$. This is a useful property of our formalism and one that is not presently shared by approaches
to this question through algebraic geometry. 
In demonstrating this property, we confine ourselves to stating the necessary results and to
giving explicit generating series for the case $k=1$ (next to minimum) and for a few values
of~$(g,n).$

\subsection{The general case}
By extending Theorem~\ref{minhur} to obtain full terms of total degree one higher than the ``minimum,''
we obtain the following result that identifies 
$\law\tau_{\al_1} \cdots \tau_{\al_n}\la_{g-1}\raw$
as a coefficient in the generating series $\La_{n,1}^g$, where
we use the notation $\La_{n,k}^g=\sF_{2g-3+2n+k} \sC\,\Xin H_n^g$ for
the terms that are $k$ higher than minimum total degree, $k\ge 0$.

\begin{tm}\label{T:k1}
For $n,g\ge 1$  and $n\ge 3, g=0$, \label{minhur2}\lremind{minhur2}
\begin{eqnarray*}
\La^g_{n,1}(\bfy ) &=& y_1\cdots y_n (-1)^{3g-3+n}\!\!\!\!
\sum_{{\be\vdash_0 2g-2+n,} \atop {l(\be)=n}}\!\!\!\!
\law\tau_{\be_1} \cdots \tau_{\be_n}\la_{g-1}\raw \be_1! \cdots \be_n!m_{\be}(\bfy) \\
&\mbox{}& +(-1)^{3g-2+n}c_g(2g-3+n)! \; y_1\cdots y_n\!\!\!\!
\sum_{k=2}^{2g-2+n} \left(\sum_{j=1}^{k-1}\frac{1}{j}\right) p_k(\bfy) h_{2g-2+n-k}(\bfy).
\end{eqnarray*}
where $\be=(\be_1, \ld  ,\be_n)$.
\end{tm}

\newcommand{\ddy}{\frac{y_1^2\partial}{\partial y_1}}
\newcommand{\pdy}{\frac{\partial}{\partial y_1}}

By extending Theorem~\ref{newpde}, we obtain a partial differential equation that is satisfied
by the generating series $\Lambda^g_{n,1}$. This is stated in the following theorem.
(Recall that $\Om^g_{n}=\Lambda^g_{n,0}$, where $\Om^g_{n}$ is used in the
proof of Theorem~\ref{newpde}.)
\begin{tm}\label{T:k1r}
For $g,n\ge1,$
$$\La^g_{n,1}(\bfy )+\frac{1}{n}
\SYM^y_{1,1} \frac{y_1^3y_2}{y_1-y_2} \frac{\partial}{\partial y_1}
\La^g_{n-1,1}(y_1,y_3,\ldots,y_n)
=\frac{1}{n} \left(T^{''}_1+ \cdots + T^{''}_4  \right)$$
where
\begin{eqnarray*}
T^{''}_1&=& 
\left(\sum_{i=1}^n \frac{y_i^2\partial}{\partial y_i}\right) \Om_{n-1}^{g}(\bfy),
\\
T^{''}_2&=&2\SYM^y_{1,1} \f{y_1^4y_2}{y_1-y_2}\pdy\Om_{n-1}^{g}(y_1,y_3, \ld  ,y_{n}),
\\
T^{''}_3&=&-\sum_{k=3}^n\,\SYM^y_{1,k-1}
\left(\ddy \Om_{k}^{0}(y_1, \ld  ,y_{k})\right)
\left(\ddy \Om_{n-k+1}^{g}(y_1,y_{k+1}, \ld  ,y_{n})\right),
\\
T^{''}_4&=& -\tfrac{1}{2}\sum_{{1\le k\le n,} \atop {1\le a\le g-1}}
\SYM^y_{1,k-1}\left( \ddy  \Om_{k}^{a}(y_1, \ld  ,y_{k})\right)
\left(\ddy\Om_{n-k+1}^{g-a}(y_1,y_{k+1}, \ld  ,y_{n})\right)
\end{eqnarray*}
where $\La_{0,1}^g = \Om^g_0 =0$ for all $g.$
\end{tm}

Note that the right hand side of the partial differential
equation for $\La^g_{n,1}$ in Theorem~\ref{T:k1r} involves only the
previously determined series $\Om^g_n=\La_{n,0}^g$.
A similar partial differential equation can be derived for the generating series
$\La_{n,k}^g$, 
given in general form in the following result.
\begin{tm}\label{genk}
For $k\ge0$,
$$\Lambda^g_{n,k}(\bfy ) + \frac{1}{n+k-1} \SYM^y_{1,1} 
\f{y_1^3y_2}{y_1-y_2}\pdy\Lambda_{n-1,k}^{g}(y_1,y_3, \ld  ,y_{n})$$
depends only upon $\Lambda^l_{j,i}$ for $0\le i<k,$ $0\le l\le g,$ $1\le j\le n.$
\end{tm}

We observe that Theorem~\ref{T:k1r} agrees with the case $k=1$, and
that equation~(\ref{dfeqzero}) agrees with the case $k=0$, in which
the right hand side is identically zero. We do not know how to
exploit the fact that the partial differential
operator applied to $\Lambda_{n-1,k}^{g}$ in Theorem~\ref{genk} is independent of $k$. 

\subsection{Explicit results for $k=1$}
\subsubsection{The genus $g=1$ case}
For the genus $g=1$ case we have the following corollary of Theorem~\ref{T:k1r}
that, together with Theorem~\ref{T:k1}, gives an explicit expression for the
generating series $\La^g_{n,1}$ for the intersection numbers
$\law\tau_{\be_1} \cdots \tau_{\be_n}\la_{g-1}\raw.$ 
\begin{co}\label{Cor:genus1}
\begin{eqnarray*}
\La^1_{n,1}(\bfy )
&=&  \frac{(-1)^{n+1}}{24}  (n-1)!  y_1\cdots y_n
\sum_{k=2}^n  \left(\sum_{j=1}^{k-1}\frac{1}{j}\right) p_k(\bfy) h_{n-k}(\bfy)
+\frac{(-1)^n}{24} n! y_1\cdots y_n h_n(\bfy) \\
&\mbox{}& +\frac{(-1)^{n-1}}{24}y_1\cdots y_n\sum_{i=2}^n
\sum_{m=i}^n \sum_{k=0}^{n-m}(i-2)!(n-i)!
(-1)^{m-i}{m\choose i} e_m(\bfy)
h_k(\bfy)h_{n-k-m}(\bfy).
\end{eqnarray*}
\end{co}

The resolutions of the generating series $24\La^1_{n,1}$, in which $g=1$,
with respect to the monomial symmetric functions $m_\theta,$ where $\theta$ is a partition,
are listed below for $1\le n\le6.$ They are obtained directly from Corollary~\ref{Cor:genus1}.
(Note that $24=c_1^{-1}$.)
$$
\begin{array}{cc|c}
g & n &24\La^1_{n,1} \\  \hline\hline
1 & 1 & -m_{2} \\
1 & 2 & m_{3\,1}+m_{2^2}  \\
1 & 3 & -m_{4\,1^2} -2m_{3\,2\,1}-2m_{2^2} \\
1 & 4 &  -2m_{5\,1^3} +3m_{4\,2\,1^2} + 4m_{3^2\,1^2} + 6m_{3\,2^2\,1} +6m_{2^4} \\
1 & 5 & 34m_{6\,1^4} +8m_{5\,2\,1^3}-12m_{4\,2^2\,1^2}-16m_{3^2\,2\,1^2}-24m_{3\,2^3\,1}-24m_{2^5}\\
1 & 6 &  -324m_{7\,1^5} -170m_{6\,2\,1^4} -112m_{5\,3\,1^4} -40m_{5\,2^2\,1^3}
               -96m_{4^2\,1^4}  \\
   &    &   +60m_{4\,2^3\,1^2}+24m_{3^3\,1^3}+80m_{3^2\,2^2\,1^2}+120m_{3\,2^4\,1}
               +120m_{2^6}.
\end{array}
$$
The intersection numbers $\law\tau_{\al_1} \cdots \tau_{\al_n}\la_{g-1}\raw$ for $g=1$ are then given
by Theorem~\ref{T:k1}.

\subsubsection{The arbitrary genus case}
The next table gives generating series $c_g^{-1} \La^g_{n,1}$ for $g=2,\dots,5$ 
and for a few values of $n.$  The series are obtained from Theorem~\ref{T:k1r}.

$$
\begin{array}{cc|c}
g & n &c_g^{-1}  \La^g_{n,1} \\  \hline\hline
2 & 1 & 37m_{4} \\
2 & 2 & -106m_{5\,1}-111m_{4\,2}-116m_{3^2}  \\
2 & 3 & 362m_{6\,1^2} +424m_{5\,2\,1}+444m_{4\,3\,2} +444m_{4\,2^2}+464m_{3^2\,2}\\ \hline
3  & 1 & -3426m_6 \\
3 & 2  & 16836m_{7\,1}+17130m_{6\,2}+17424m_{5\,3}+17424m_{4^2} \\ \hline
4 & 1 & 61164m_8 \\
4 & 2 & -4249232m_{9\,1} -4278148m_{8\,2} -4307064m_{7\,3} -4311180m_{6\,4} 
              -4315296m_{5^2} \\ \hline
 5 & 1 & -180519696m_{10} \\
 5 & 2  & 1619765280 m_{11\, 1}+1624677264 m_{10\, 2}+1629589248 m_{9\, 3} \\
   &    &   +1630276704 m_{8\, 4}+1630964160 m_{7\, 5}+1630964160 m_{6^2}.
\end{array}
$$

} 


\begin{thebibliography}{[FabP2]}
\bibitem[ELSV1]{elsv1} T. Ekedahl, S. Lando, M. Shapiro, and A. Vainshtein,
{\em On Hurwitz numbers and Hodge integrals}, C. R. Acad.\ Sci.\ 
Paris S\'{e}r. I Math.\ {\bf 328} (1999), 1175--1180.  \lremind{elsv1}  \cited

\bibitem[ELSV2]{elsv2} T. Ekedahl, S. Lando, M. Shapiro, and A. Vainshtein,
{\em Hurwitz numbers and intersections on moduli spaces of curves},
Invent.\ Math.\ {\bf 146} (2001), 297--327. 
 \cited

\bibitem[Fab]{fabconj}
C. Faber, {\em A conjectural description
of the tautological ring of the moduli space of curves}, in 
{\em Moduli of Curves and Abelian Varieties}, 109--129, Aspects Math.,
{\bf E33}, Vieweg, Braunschweig, 1999.   \lremind{faber} \cited

\bibitem[FabP1]{fabplambdag} C. Faber and R. Pandharipande,
{\em Hodge integrals, partition matrices, and the $\la_g$ conjecture},
Ann.\  Math.\ (2) {\bf 157} (2003), no.\ 1, 97--124. \cited

\bibitem[FabP2]{fabpjems} C. Faber and R. Pandharipande,
{\em Relative maps and tautological classes}, J. Eur.\ Math.\ Soc.\ {\bf 7}
(2005), no.\ 1, 13--49.  \cited 

\bibitem[GeP]{gp} E. Getzler and R. Pandharipande,
 {\em Virasoro constraints and the Chern classes of the Hodge bundle},
 Nuclear Phys.\ B {\bf 530} (1998), 701--714.  \cited


 \bibitem[Giv]{giv}
A. Givental, {\em Gromov-Witten invariants and 
quantization of quadratic hamiltonians}, Mosc.\ Math.\ J. {\bf 1} (2001), no.\ 4,
551--568, 645. 
\lremind{math.AG/0108100}  \cited

\bibitem[GJ]{gjce} I. P. Goulden and D. M. Jackson,
Combinatorial Enumeration, Wiley, New York, 1983 (Reprinted by Dover, 2004).

\bibitem[GJ1]{gj0} I. P. Goulden and D. M. Jackson, 
{\em Transitive factorizations into transpositions and holomorphic
mappings on the sphere},
Proc.\  Amer.\  Math.\  Soc.\
{\bf 125} (1997), 51--60.

\bibitem[GJ2]{gjconj}
I. P. Goulden and D. M. Jackson, {\em The
number of ramified coverings of the sphere by the double
torus, and a general form for higher genera}, J. Combin.\ Theory
A {\bf 88} (1999) 259--275. \cited

\bibitem[GJVai]{gjvai}
I.P. Goulden, D.M. Jackson and A. Vainshtein,
{\em The number of ramified coverings of the sphere by the torus and
surfaces of higher genera}, Ann.\  Combinatorics {\bf 4} (2000),
27--46. \cited


\bibitem[GJV1]{gjv1} I. P. Goulden, D. M. Jackson, R. Vakil,
 {\em The Gromov-Witten potential of a point, Hurwitz numbers, and Hodge
 integrals}, Proc.\ London Math.\ Soc.\ (3) {\bf 83} (2001), 563--581. 
\cited

\bibitem[GJV2]{double} I. P. Goulden, D. M. Jackson, R. Vakil,
 {\em Towards the geometry of double Hurwitz numbers},
 Adv.\ Math.\ {\bf 198} (2005), 43--92.  \cited

\bibitem[GJV3]{gjvfaber} I. P. Goulden, D. M. Jackson, R. Vakil,
{\em The moduli space of curves, double Hurwitz numbers,
and Faber's intersection number conjecture}, in preparation. \cited

\bibitem[GP]{vl} T. Graber and R. Pandharipande,
{\em Localization of virtual classes}, Invent.\ Math.\ {\bf 135} (1999), no.\ 2, 487--518.
\cited 

 \bibitem[GV1]{socle} T. Graber and R. Vakil, {\em On the tautological
 ring of $\cmbar_{g,n}$}, in {\em Proceedings of the Seventh 
 G\"okova Geometry-Topology Conference 2000}, International Press, 2000. 
\cited

\bibitem[GV2]{gvelsv} T. Graber and R. Vakil, {\em Hodge integrals,
Hurwitz numbers, and virtual localization}, Compositio Math.\  {\bf 135} (1) 
(January 2003), 25--36.  \cited

 \bibitem[GV3]{thmstar} T. Graber and R. Vakil, {\em Relative
 virtual localization and vanishing of tautological classes on
 moduli spaces of curves}, Duke Math.\ J. {\bf 130} (2005), no.\ 1, 1--37.
\cited

\bibitem[H]{hur} A. Hurwitz, 
{\em \"{U}ber Riemann'sche Fl\"{a}chen mit gegebenen Verzweigungspunkten},
Matematische Annalen,
{\bf39} (1891), 1--60.

\bibitem[KaL]{kal} M. E. Kazarian and S. K. Lando, {\em An algebro-geometric
proof of Witten's conjecture}, Max-Planck Institute preprint MPIM 2005--55 (2005),
http://www.mpim-bonn.mpg.de/preprints. \cited

 \bibitem[KiL]{kil} 
Y.-S. Kim and K. Liu, {\em A simple proof of Witten conjecture through
localization}, preprint 2005, math.AG/0508384.\cited

 \bibitem[LP]{leep} 
Y.-P. Lee and R. Pandharipande, {\em Frobenius manifolds, 
Gromov-Witten theory, and Virasoro constraints}, book in preparation.
\cited

\bibitem[LLZ1]{llz1} C.-C. M. Liu, K. Liu, and J. Zhou, {\em
    A proof of a conjecture of Mari\~no-Vafa on Hodge integrals}, 
  J. Diff.\ Geom.\ {\bf 65} (2003), no.\ 2, 289--340. \cited

\bibitem[LLZ2]{llz} C.-C. M. Liu, K. Liu, and J. Zhou, {\em
    Mari\~no-Vafa formula and Hodge integral identities}, J.
  Algebraic Geom.\ {\bf 15} (2006), 379--398. \cited

\bibitem[MVaf]{mv}
M.Mari\~no and C.Vafa,
{\em Framed knots at large~N},
in {\em Orbifolds in Mathematics and Physics (Madison, WI 2001)},
185--204,
Contemp. \ Math., 310, Amer.\ Math.\  Soc., Providence, RI, 2002.\cited

\bibitem[OP]{op} A. Okounkov and R. Pandharipande, {\em Virasoro constraints
for target curves}, preprint 2003, math.AG/0308097. \cited

\bibitem[P]{icm} R. Pandharipande, {\em Three questions in Gromov-Witten theory},
in {\em Proceedings of the International Congress of Mathematicians, Vol.\ II (Beijing, 2002)},
503--512, Higher Ed.\ Press, Beijing, 2002. \cited

\bibitem[SZ]{sz}  S. Shadrin and D. Zvonkine, {\em Changes
of variables in ELSV-type formulas}, preprint 2006, math.AG/0602457.\cited


\bibitem[V]{cime} R. Vakil, {\em The moduli space of curves and
     Gromov-Witten theory}, preprint 2006, math.AG/0602347v2,
     submitted for publication.  \cited



\bibitem[X]{x}
G. Xin, {\em A fast algorithm for MacMahon's partition analysis},
Elec.\ J.\ Comb.\ {\bf 11} (2004), R58.\cited 

 \end{thebibliography}
 \end{document}